\documentclass[a4paper,12pt]{article}
\usepackage{amssymb,amsmath,amsthm,latexsym}
\usepackage{amsfonts}
	\usepackage{amsfonts}
\usepackage{graphicx}
\usepackage[pdftex,bookmarks,colorlinks=false]{hyperref}
\usepackage{tikz}
\usepackage{verbatim}
\usepackage{caption}
\usepackage{subcaption}

\newtheorem{theorem}{Theorem}[section]

\newtheorem{definition}[theorem]{Definition}

\newtheorem{lemma} [theorem]{Lemma}

\newtheorem{remark}[theorem]{Remark}

\title{\bf ON WEAKLY UNIFORM INTEGER ADDITIVE SET-INDEXERS OF GRAPHS}
\author{{\bf K A Germina\footnote{Department of Mathematics, School of Mathematical \& Physical Sciences, Central University of Kerala, Kasaragod, email:{\em srgerminaka@gmail.com}}}     and    {\bf N K Sudev \footnote{Department of Mathematics, Vidya Academy of Science \& Technology, Thalakkottukara, Thrissur - 680501, email: {\em sudevnk@gmail.com}}}}
\date{}

\begin{document}
\maketitle

\begin{abstract} 
We have the notion of set-indexers, integer additive set-indexers and $k$-uniform integer additive set-indexers of graphs. In this paper, we initiate a study of the graphs which admit $k$-uniform integer additive set-indexers and introduce the notion of weakly uniform integer additive set-indexers and arbitrarily uniform integer additive set-indexers and provide some useful results on these types of set-indexers.
\end{abstract} 



\section{Introduction}. 
For all  terms and definitions, not defined specifically in this paper, we refer to \cite{FH}. Unless mentioned otherwise, all graphs considered here are simple, finite and have no isolated vertices.

The set-labeling of a graph is defined in \cite{A1}as follows. Let $G$ be a $(p,q)$-graph and let $X$, $Y$,$Z$ be non-empty sets. Then the functions $f:V(G)\to 2^X$, $f:E(G)\to 2^Y$ and $f:V(G)\cup E(G)\to 2^Z$ are called the {\em  set-assignments} of vertices, edges and elements of $G$ respectively. By a set-assignment of a graph, we mean any one of them.  A set-assignment is called a {\em set-labeling} if it is injective.

\begin{definition}\label{D0}{\rm
\cite{A1} For a $(p,q)$- graph $G=(V,E)$ and a non-empty set $X$ of cardinality $n$, a {\em set-indexer} of $G$ is defined as an injective set-valued function $f:V(G) \rightarrow2^{X}$ such that the function $f^{\oplus}:E(G)\rightarrow2^{X}-\{\emptyset\}$ defined by $f^{\oplus}(uv) = f(u ){\oplus} f(v)$ for every $uv{\in} E(G)$ is also injective, where $2^{X}$ is the set of all subsets of $X$ and $\oplus$ is the symmetric difference of sets.}
\end{definition}

\begin{theorem}
\cite{A1} Every graph has a set-indexer.
\end{theorem}

\begin{definition}\label{D1}{\rm
Let $N_0$ denote the set of all non-negative integers. For all $A, B \subseteq N_0$, the sum of these sets is denoted by  $A+B$ and is defined by $A + B = \{a+b: a \in A, b \in B\}$}.
\end{definition}

\begin{definition}\label{D2}{\rm
\cite{GA} An {\em integer additive set-indexer} (IASI, in short) is defined as an injective function $f:V(G)\rightarrow 2^{N_0}$ such that the induced function $g_f:E(G) \rightarrow 2^{N_0}$ defined by $g_f (uv) = f(u)+ f(v)$ is also injective}.
\end{definition}

\begin{definition}\label{D3}{\rm
The cardinality of the labeling set of an element (vertex or edge) of a graph $G$ is called the {\em set-indexing number} of that element.}
\end{definition}

\begin{definition}\label{DU}
\cite{GA} An IASI is said to be {\em $k$-uniform} if $|g_f(e)| = k$ for all $e\in E(G)$. That is, a connected graph $G$ is said to have a $k$-uniform IASI if all of its edges have the same set-indexing number $k$. In particular, we say that a graph $G$ has an {\em arbitrarily $k$-uniform IASI} if $G$ has a $k$-uniform IASI  for every positive integer $k$.
\end{definition}

In \cite{GA}, a study of a special type of 1-uniform IASI {\em sum square graphs} and in\cite{TMKA}, a characterisation of 2-uniform IASI were established. The most important result given in \cite{TMKA} is the following.

\begin{theorem}\label{T0}
\cite {TMKA} A graph $G$ has a 2-uniform IASI if and only if it is a bipartite graph.
\end{theorem}

The motivation for the study presented in this paper is from the results in \cite{GA} and \cite{TMKA}. As a generalisation of 1-uniform and 2-uniform IASIs, in section 2, we introduce the notion of weakly $k$-uniform IASI and arbitrarily $k$-uniform IASI and prove the conditions for a connected graph to admit a weakly or arbitrarily $k$-uniform IASI, where $k$ is a positive integer. 
 
\section{Uniform Integer Additive Set-Indexers}

During this study our first attempt is to analyse the set-indexing number of an edge of a graph $G$ and to study its relation with the set-indexing number of the end vertices of that edge . As a result,  we establish the following results regarding the cardinality of $A+B$ defined in Definition \ref{D1}. If either $A$ or $B$ is countably infinite, then clearly $A+B$ is also countably infinite and hence the study of the cardinality of $A+B$ trivial. Hence we restrict our discussion  on finite sets $A$ and $B$. We denote the cardinality of a set $A$ by $|A|$. 

\begin{lemma}\label{T1}
Let $A,B\subseteq N_0$. Then $max(|A|,|B|)\le |A+B|\le|A||B|$.
\end{lemma}

{\bf Proof:} Let $A,B \subseteq N_0$. Also let $|A|=m,|B|=n$. Assume either $A$ or $B$ is a singleton set. Without loss of generality, let $A$ be a singleton set, say $A=\{q\}$ where $q$ is a non-negative integer. That is, $m=1$. Therefore, $A+B=\{b_1+q, b_2+q,.. \cdots b_n+q\}$. Hence $|A+B|=n$. Similarly, if $|B|=1$, then $|A+B|=m$. Hence, $max(|A|,|B|)\le |A+B|$. If  neither $A$ nor $B$ is singleton, then define $h:A\times B \to A+B$ defined by, $h(a,b)=a+b ; a\in A,b\in B$. Clearly, $h$ is injective and hence $|A+B| \le |A||B|$. The equality holds when the function $h: A \times B \to A+B$ is a bijection. From the above two cases, we have $max(|A|,|B|)\le |A+B|\le|A||B|$.

\begin{remark}\label{R1}{\rm
Due to Lemma \ref{T1}, given an integer additive set-indexer $f$ of a graph $G$,  $max(|f(u)|, |f(v)|) \le |g_f(uv)|= |f(u)+f(v)| \le |f(u)| |f(v)|$, where $u,v \in V(G)$.}
\end{remark}

The characteristics of the IASIs with the property $|A+B|=max(|f(u)|,\\|f(v)|)$ is of special interest. Hence we have the following definition.

\begin{definition}{\rm
An IASI $f$ of a graph $G$ is called a {\em weak IASI} if $|g_f(uv)| = max(|f(u)|,|f(v)|)$ for all $u, v \in V(G)$. } 
\end{definition}

It is  observed that if $|g_f(uv)| = max(|f(u)|,|f(v)|) \forall u,v\in V(G)$ holds for every edge $uv$ of $G$, it is necessary that one of its end vertices is labeled by a singleton set. Hence we define a {\em weakly $k$-uniform IASI} as follows. 

\begin{definition}{\rm
A $k$-uniform IASI which assigns only singleton sets and $k$- element sets to the vertices of a given graph $G$ is called a {\em weakly $k$-uniform IASI}.}
\end{definition}

\begin{theorem}\label{TT}
Every tree admits a weakly $k$-uniform IASI. 
\end{theorem}

{\bf Proof:} All the sets under consideration are subsets of $N_0$. Let $T$ be the given tree and $k$ be an arbitrary positive integer. If $k=1$, assign distinct singleton sets to every vertex of $T$. Then clearly, every edge of $T$ also has singleton sets as labeling sets. Therefore, $T$ is 1-uniform. 
\\Now let $k>1$. 

{\em Step-1:} Let $u_1,u_2, \cdots u_r$ be the pendant vertices of $T$. Assign the singleton sets  $\{i\}$ to the pendant vertices  $u_i, 1\le i\le r$. Let $\{v_1,v_2,\cdots v_s\}$ be the set of internal vertices of $T$ which are adjacent to some of $u_i,1\le i\le s$. Assign the $k$-element set $\{j,j+1,j+2,\cdots j+(k-1)\}$ to each vertex $v_j,1\le j\le s$. Then, by Theorem \ref{T1}, each edge $u_iv_j$ has the set-indexing number $k$. Note that no $u_i$ can be adjacent to two internal vertices which are adjacent to each other, since otherwise, we have a contradiction to the acyclic property of $T$. 

{\em Step-2:} Let $\{w_1,w_2,\cdots w_l\}$ be the set of internal vertices which are adjacent to some of $v_j,1\le j\le l$,   where no $w_n$ can be adjacent to any pendant vertices. Assign singleton sets $\{r+n\}$ to each $w_n, 1\le n\le l$ so that each edge $v_jw_n$ has the set-indexing number $k$. 

The succeeding steps, the vertices that are adjacent to the just preceding labeled vertices alternately by distinct $k$-element sets other than the already labelled  sets and distinct singleton sets, for a finite number of times, we get a weakly $k$-uniform IASI to the tree $T$.

\begin{remark}{\rm
A path $P_n$ can be considered as a tree which has two pendant vertices. Hence by Theorem \ref{TT}, $P_n$ admits an arbitrary $k$-uniform IASI.}
\end{remark}

\begin{theorem}\label{TCE}
Every even cycle $C_{2n}$ has a weakly $k$-uniform IASI.
\end{theorem}

{\bf Proof:} Let $C_{2n}=v_1v_2v_3\cdots v_{2n}v_1$ be a cycle of even length. Assign distinct singleton subsets of $N_0$ to the vertices $v_1,v_3,\dots v_{2n-1}$ and assign distinct $k$-element subsets of $N_0$ to the vertices $v_2, v_4,\cdots v_{2n}$ of $C_{2n}$. Then by Theorem \ref{T1}, each edge $v_iv_j$ has the set-indexing number $k$. Hence an even cycle has an arbitrary $k$-uniform IASI.

\begin{theorem}\label{TCO}
No odd cycle $C_n$ has a weakly $k$-uniform IASI.
\end{theorem}

{\bf Proof:} By the definition of weakly  $k$-uniform IASI, for each edge $v_iv_{i+1}$ in $C_n$, one of $v_i$ (or $v_{i+1}$) should necessarily be  labeled by a singleton set and $v_{i+1}$ (or $v_i$) with a $k$-element set. Since $n$ is odd, $n=2m+1$ for some non-negative integer $m$. Let $C_{n}=v_1v_2v_3\cdots v_{2m+1}v_1$ be a cycle of odd length. Assign distinct singleton subsets of $N_0$ to the vertices $v_1,v_3,\dots v_{2m-1}$ and assign distinct $k$-element subsets of $N_0$ to the vertices $v_2, v_4,\cdots v_{2m}$ of $C_{n}$. Then by Theorem \ref{T1}, each edge $v_iv_j$ where $i,j\le 2n$, has the set-indexing number $k$. Now if we assign a singleton set to the vertex $v_{2m+1}$, then the edge $v_{2m+1}v_0$ has the set-indexing number 1. If we assign a $k$-element set to $v_{2m+1}$, then the edge $v_{2m}v_{2m+1}$ has the set-indexing number greater than $k$. Hence no odd cycle admits a weakly $k$-uniform IASI.

Invoking the above results, the existence of a weakly $k$-uniform IASI implies the bipartiteness of the given graph.

\begin{theorem}\label{TBG}
Let $k$ be any positive integer. A bipartite graph $G$ admits a weakly $k$-uniform IASI.
\end{theorem}

{\bf Proof:} Let $G=(X_1,X_2)$ be bipartite graph. Then assign distinct singleton sets to the vertices in $X_1$ and distinct $k$-element sets to $X_2$ so that every edge of $G$ has the set-indexing number $k$. Therefore, $G$ admits a weakly $k$-uniform IASI.

\begin{theorem}\label{TBW}
For any positive integer $k>1$, a graph $G$ admits a weakly $k$-uniform IASI if and only if $G$ is bipartite.
\end{theorem}

{\bf Proof:} Proof of the necessary part follows from \ref{TBG}. Conversely, assume that $G$ admits a weakly $k$-uniform IASI. That is, every edge of $G$ has the set-indexing number $k$. Since $G$ admits a weakly $k$-uniform IASI, one end vertex of all edges of $G$ should necessarily be labeled by distinct singleton sets and the other end vertex by distinct $k$-element sets. Let $X_1=\{u_1,u_2,u_3, \cdots, u_r\}$ be the set of vertices which have been labeled by distinct singleton sets and let $X_2=\{v_1,v_2,v_3,\cdots v_s\}$ be the set of vertices which have been labeled by $k$-element sets. Since $k>1$, both $X_1$ and $X_2$ are non-empty and no two vertices of $X_1$ (and $X_2$) can be adjacent to each other. That is, $(X_1,X_2)$ is a bipartition of the vertex set of $G$.

Next we establish the existence of $k$-uniform IASIs to bipartite graphs in the following theorem.

\begin{theorem}\label{TB}
A bipartite graph $G$ admits a $k$-uniform IASI.
\end{theorem}

{\bf Proof:} Let $G$ be a bipartite graph. Let $(X_1, X_2)$ be the bipartition of $V(G)$. Let $X=\{v_1,v_2,v_3,\cdots, v_l\}$ and $Y=\{v_{l+1},v_{l+2},v_{l+3},\cdots, v_r\}$.  Let $m,n,d$ be positive integers with $m,n>1$ and $d\ge 1$. Now assign the set $\{i,i+d,i+2d, \cdots, i+(m-1)d\}$ to the vertex $v_i,1\le i\le l$ in $X_1$. Let $v_{l+j},j\ge 1$, be a vertex in $X_2$ which is adjacent to a vertex $v_i$. Now assign the set $\{l+j, l+j+d, l+j+2d,\cdots, l+j+(n-1)d\}$ to the vertex $v_{l+j}$. Now the edge $v_iv_{l+j}$ has the set-indexer $\{i+j+l, i+j+l+d, i+j+l+2d, \cdots i+j+l+(n-1)d, i+j+l+d, i+d+j+l+d, i+d+j+l+2d, \cdots i+d+j+l+(n-1)d, \cdots, i+(m-1)d+j+l, i+(m-1)d+j+l+d, \cdots, i+(m-1)d+j+d+(n-1)d\}$. That is,  the set-indexer of $v_iv_{l+j}$ is $\{i+j+l,i+j+l+d,i+j+l+2d, \cdots, i+j+l+(m-1)d, i+j+l+md, i+j+l+(m+1)d, i+j+l+(m+2)d, \cdots, i+j+l+(m+n-2)d\}$.
Therefore, the set-indexing number of $u_iv_{l+j}$ is $m+n-1$. Since $v_i$ and $v_{l+j}$ are arbitrary elements of $X_1$ and $X_2$ respectively, the above argument can be extended to all edges in $G$. That is, all edges of $G$ can be assigned by a $k$-element set, where $k=m+n-1$. Therefore, $G$ admits a $k$-uniform IASI.

\begin{theorem}\label{TCW}
A complete graph $K_n, n>2$ does not admit a weakly $k$-uniform IASI for any positive integer $k>1$.
\end{theorem}

{\bf Proof:} Let $k>1$ be a positive integer. If possible, let $K_n$ admits a weakly $k$-uniform IASI. Then its vertices are assigned by either singleton sets or $k$-element sets. As explained in the proof of Theorem \ref{TBW}, we get a bipartition $(X_1,X_2)$ of the vertex set of $K_n$. This is a contradiction to the fact that no complete graph other than $K_2$ is bipartite. Therefore, $K_n, n>2$ does not admit a weakly $k$-uniform IASI.

It is interesting to check whether implications of the existence of a weakly $k$-uniform IASI to the existence of an arbitrary $k$- uniform IASI. For any positive integer $k$, if there exists a weakly $k$-uniform IASI for a graph $G$, then each of its vertices is labeled by either a singleton set or by a $k$-element set and all edges of $G$ have the set-indexing number $k$. Therefore, $G$ has a $k$-uniform IASI. Since $k$ is an arbitrary positive integer, we arrive at the following theorem.

\begin{theorem}\label{TWA}
A graph which has a weakly $k$-uniform IASI admits an arbitrarily $k$-uniform IASI.
\end{theorem}

Invoking Theorem \ref{TWA}, every path, even cycle and tree have an arbitrarily $k$-uniform IASI. Hence we can write the following theorem. 

\begin{theorem}\label{TBA}
Every bipartite graph admits an arbitrarily $k$-uniform IASI, where $k$ is a positive integer.
\end{theorem}

We have already established the existence of $k$-uniform IASIs to bipartite graphs for any positive integer $k$. Hence we shall look in to  the existence of such uniform IASIs for other connected graphs. 

It can be observed that the completely bipartite graph $K_{m,n}$ admits a $k$-uniform IASI for any positive integer $k$. For any given positive integer $k\ge 1$, the existence of $k$-uniform IASI is established for a bipartite graph. Now recall the Theorem \ref{T0} that a connected graph $G$ admits a 2-uniform IASI if and only if it is bipartite. Hence no non-bipartite graphs can have a 2-uniform IASI. Then, what are the admissible values of $k\ne 2$, for a non-bipartite graph to have a $k$-uniform IASI? The following results settle this problem.

\begin{theorem}\label{TS}
Let $k$ be a positive integer. If $G$ is a graph which admits a $k$-uniform IASI, then any subgraph $H$ of $G$ also admits the same $k$-uniform IASI.
\end{theorem}

{\bf Proof:} Let $G$ be a graph which admits a $k$-uniform IASI and $H$ be a subgraph of $G$. Let $f^\ast$ be the restriction of $f$ to $V(H)$. Then $g_{f^\ast}$ is the corresponding restriction of $g_f$ to $E(H)$. Then clearly, $f^\ast$ is a set-indexer on H. This set-indexer may be called the {\em induced set-indexer} on $H$. Since $g_f(e) = k$ for all $e \in E(G)$, we have $g_{f^\ast}(e) =k$ for all $e \in E(H)$. Hence $H$ admits $k$-uniform IASI.

\begin{remark}\label{Tsup}{\rm
As the contrapositive statement of the above theorem, we observe that if a graph $G$ does not have a $k$-uniform IASI, then any supergraph of $G$ does not admit a $k$-uniform IASI. }
\end{remark}

\begin{theorem}\label{TC}
A complete graph $K_n$ admits a $k$-uniform IASI for positive odd integer $k$.
\end{theorem}

{\bf Proof:} Let $\{v_1,v_2,v_3,\cdots,v_n\}$ be the vertex set of the complete graph $K_n$. Let $d$ and $m$ be a positive integer. Assign the set $\{i,i+d,i+2d,\cdots,i+(m-1)d\}, i\le i\le n$ to each vertex $v_i$. Then an edge $v_iv_j$ of $K_n$ has the indexing set $\{i+j,i+j+d,i+j+2d,\cdots,i+j+2(m-1)d\}$.That is, the set-indexing number of $v_iv_j$ is $2m-1$. Since $v_i$ and $v_j$ are arbitrary vertices of $K_n$, all edges of it has the set-indexing number $2m-1$. Then $K_n$ admits a $k$-uniform IASI, where $k=2m-1$.

\begin{remark}\label{R2}{\rm
Being a subgraph of a complete graph $K_n$,  any connected graph $G$ admits a $k$-uniform IASI for all positive odd integers $k$.}
\end{remark}

Hence the interest of studying  the graphs which admit $k$-uniform IASI may be related to an even integer $k>2$. Hence what are these graphs which admit $k$-uniform IASI for even integer $k$, other than bipartite graphs? The following theorem settles this problem.

\begin{theorem}\label{TCNB}
Let $G$ be a connected non-bipartite graph. Then $G$ admits a $k$-uniform IASI if and only if $k$ is a positive odd integer.
\end{theorem}

{\bf Proof:} Necessary part of the result follows from Remark \ref{R2}. Let $m,n$ be two positive integers. Let $G=(X_1,X_2)$ be a bipartite graph. Assign $m$-element sets to the vertices of $X_1$ and $n$-element sets to the vertices of $X_2$ as in Theorem \ref{TB} so that every edge of $G$ has the set-indexing number $k=m+n-1$. Draw an edge between two vertices of $X_1$ (or $X_2$) that gives a minimal connected non-bipartite graph with given number of vertices. Without loss of generality, let $e=v_iv_j$ where $v_i,v_j \in X_1$ (or $e=v_{l+i}v_{l+j};v_{l+i},v_{l+j}\in X_2$. Let $G_1=G\cup \{e\}$. If possible, assume that $G_1$ admits a $k$-uniform IASI where $k>2$ is an even integer. Then, the set-indexing number of $e$ is $2m-1$ (or $2n-1$). Since $G_1$ admits a $k$-uniform IASI, $k>2$, we have $m+n-1=2m-1=2n-1=k$, which is possible on when $m=n$, a contradiction. Hence $k$ can not be even. Hence, for a positive even integer $k$, by Remark \ref{Tsup}, no supergraph of $G_1$ (including $K_n$) admits a $k$-uniform IASI for a positive integer $k$.

Hence, a graph $G$ admits a $k$-uniform IASI, for a positive even integer $k>2$ if and only if $G$ is bipartite. More generally, 

\begin{theorem}
A graph $G$ admits a $k$-uniform IASI, for a positive integer $k$, if and only if $k$ is odd or $G$ is bipartite.
\end{theorem}

\end{document}